\documentclass[11pt]{amsart}
\usepackage{fullpage}
\usepackage{amsmath,amssymb,amscd, pb-diagram, amssymb}

\def\symbdown#1{\Big\downarrow\rlap{$\vcenter{\hbox{$\scriptstyle
#1$}}$}}

\newtheorem{thm}{Theorem}
\numberwithin{thm}{section}

\newtheorem{lma}[thm]{Lemma}
\newtheorem{prop}[thm]{Proposition}

\newtheorem{defn}[thm]{Definition}

\begin{document}

\begin{abstract}
We show that the family of nest algebras with $r$ non-zero nest
projections is stable, in the sense that an approximate containment of
one such
algebra within another is close to an exact containment. We use this
result
to give a local characterisation of
limits formed from this family. We then consider quite general regular
limit algebras and
characterise these algebras using a local condition which reflects the
assumed regularity of the system. 
\end{abstract}

\bibliographystyle{amsplain}
\title{Local characterisation of approximately finite operator algebras}
\author{P.A. Haworth}
\date{15 December 1999}
\maketitle
\noindent
\section{Introduction}
The approximately finite (AF) $C^{*}$-algebras are completely
characterised among separable $C^{*}$-algebras by the local description
that any finite family of
elements almost lies in a finite dimensional $C^{*}$-subalgebra. This
was
proved first by Glimm~\cite{glimm} for the unital UHF algebras and later
extended
by  Dixmier and Bratteli; for details
consult~\cite{bratteli},~\cite{dixmier}. More recently,
Heffernan~\cite{hef} generalised these results to non-selfadjoint
contexts and
showed that  the  uniformly $T_{2}$-algebras, (limits of nest
algebras with self-adjoint part $M_{n_{1}}\oplus M_{n_{2}}$) admit a
similar characterisation. The first half of this paper,
will be concerned with extending this intrinsic characterisation to
cover the AF nest algebras of bounded diameter, or, in keeping with
the above  terminology, the uniformly $T_{r}$-algebras
for arbitrary fixed $r$. This paper then addresses the
characterisation problem posed in~\cite{scp2}.
For $r$ exceeding 2, star extendible embeddings
between $T_{r}$-algebras need not be decomposable into multiplicity
one embeddings, that is, they need not be regular. Accordingly we need
quite different methods from those of~\cite{hef} .

In $C^{*}$-algebra theory, one can use functional calculus techniques
to show that the family of
finite dimensional $C^{*}$-algebras is a stable one, in the sense
that an approximate inclusion of one finite dimensional $C^{*}$-algebra
in another can
be perturbed to a nearby exact inclusion, see~\cite{scp}. This property
is more
elusive for general non-selfadjoint families, yet provides a
sufficient condition for a Glimm style characterisation of limits
formed from such a family. The central result then of the first part
is to show that the $T_{r}$-algebras form a stable family. Solving the
question of stability is typical of the perturbational problems we
have to resolve. The philosophy is that: if a property is approximately
true of something, is it close to something similar for which the
property holds exactly. This theme is well developed in
$C^{*}$-algebras (see~\cite{ch} ~\cite{lor},for example), but less so in
the
framework of non-selfadjoint operator algebra.

Later, we shall focus attention on those
algebras arising as the dense union of a chain of digraph algebras
each regularly embedded in the next, a redundant assumption for
$C^{*}$-algebras and  $T_{2}$-algebras. The regular star
extendible embeddings are the most tractable mappings between digraph
algebras, essentially carrying matrix units over to sums of matrix
units. In
this  setting we
provide a new local description which reflects the assumed regularity
of the system. More importantly however, we will be able to dispense
with the bounded diameter constraint imposed above and
considerably widen the class of algebra we characterise.

The techniques used in the regular setting are of necessity quite
different from those used in the general star extendible case. We
extend the notion of a normalising partial isometry to that of an
approximately normalising, approximate partial isometry (see
Definitions 2.2 and 3.1)  and show that such an element is close to
an exactly normalising partial isometry, with the closeness depending
not on
the containing algebra, but on how well the element normalises the
masa. It is effectively this lack of dependence on the containing
algebra that will allow us to unbound the diameter of the building
block algebras. The proof of this result requires  an application
of Arveson's distance formula. 

Throughout the paper, all algebras will be assumed separable without
further mention. Uniform limits of
digraph algebras will be
taken with respect to star extendible algebra embeddings, with no
further assumption, until section 3, where they shall be taken to be
regular.  A symbol of the form $\delta(\epsilon)$ will be
taken to
denote a positive function of $\epsilon$ with the property that
$\delta(\epsilon)\rightarrow0$ as $\epsilon \rightarrow0.$ We adopt
this convention to prevent unnecessary notation in proofs.

\section{Stability for nest algebras}

We shall say that a family of algebras $ \mathcal{F}$, is
stable  if, given $\epsilon >0$ and $\mathcal{A}_{1}\in \mathcal{F}$
there
exists $\delta >0$ such that whenever $\mathcal{A}_{2}\in \mathcal{F} $
and $\phi
_{1}:\mathcal{A}_{1}\rightarrow M_{n}$ , $\phi
_{2}:\mathcal{A}_{2}\rightarrow M_{n}$  are
star extendible embeddings with $\phi _{1}(\mathcal{A}_{1})\subseteq
_\delta \phi _{2}(\mathcal{A}_{2}),$ then there exists a star extendible
algebra
injection $\psi :\mathcal{A}_{1}\rightarrow M_{n}$ with $\left\| \phi
_{1}-\psi
\right\| <\epsilon $ and $\psi (\mathcal{A}_{1})\subseteq
\phi_{2}(\mathcal{A}_{2})$. 

As alluded to
in the introduction, uniform limits formed from algebras within a
stable family may be locally characterised. More formally, let $
\mathcal{F}$ be a stable family of finite dimensional operator
algebras and let $ \mathcal{A}$ be a Banach subalgebra of a
$C^{*}$-algebra. Then the following are equivalent;
\begin{enumerate}

\item There exists a chain $ \mathcal{A}_{1}\subseteq
  \mathcal{A}_{2}\subseteq \ldots$ of subalgebras of $\mathcal{A}$
  with  $\mathcal{A}=cl{(\bigcup_{i}\mathcal{A}_{i})}$ and with  each
  $\mathcal{A}_{i}$ star extendibly isomorphic to an algebra in $
  \mathcal{F}$.

\item For each $ \epsilon>0$ and finite subset $ \mathcal{S}\subseteq
  \mathcal{A}$ there exists a pair $(\mathcal{B},\phi)$ of a finite
  dimensional operator algebra $ \mathcal{B}\in\mathcal{F}$ and a star
  extendible injection $\phi: \mathcal{B}\rightarrow \mathcal{A}$ with
  $dist( \mathcal{S},\phi(\mathcal{B}))<\epsilon$.
 
\end{enumerate}
The equivalence of the above conditions is routine given the stability
condition on the family $\mathcal{F}$. To check a family for stability
however,
is a non trivial matter and the main achievement of Glimm's
fundamental paper~\cite{glimm} was to show that the family of matrix
algebras is stable. 

The algebras we wish to characterise are the uniformly
   $T_{r}$-algebras, for arbitrary fixed $r$.   These algebras
   manifest themselves as uniform limits formed from algebras in the
family
$\mathcal{N}_{r}$ of all nest algebras with $r$ non-zero nest
projections. The
local characterisation will follow immediately from the next theorem,
the proof of which will occupy us for the remainder of this section.
\begin{thm}
$\mathcal{N}_{r}$ is a stable family.
\end{thm}
The following definition, of an
approximate partial isometry, is in the preturbational spirit and sets
up the discussion to follow.

\begin{defn}
Let $v$ be an operator on a Hilbert space $\mathcal{H}$. Then $v$ is
said to be an $%
\epsilon $-approximate partial isometry if $\left\| v^{\ast }v-(v^{\ast
}v)^{2}\right\| \leq \epsilon .$ Operators $v,w\in B(\mathcal{H})$ are
said to be $%
\epsilon $-approximately orthogonal, if $\left\| v^{\ast }w\right\| \leq
\epsilon \;$and $\left\| vw^{\ast
}\right\| \leq \epsilon .$
\end{defn}

\begin{lma}
Let $b$ be a self adjoint element of a finite dimensional
$C^{*}$-algebra  with $\left\| b^{2}-b\right\| \leq \delta ,\;\delta
<\frac{1}{4}.$
Then there exists a projection $p\in C^{\ast }(b)$ with $\left\|
p-b\right\|
<2\delta .$
\end{lma}
This result is well known and easily proved by appealing to
the functional calculus.

\begin{lma}
Let $v$ be an operator on a finite dimensional Hilbert space
$\mathcal{H}$. Given $%
\delta >0\;$ there exists $\;\epsilon >0$ such that if $v$ is an
$\epsilon $-approximate  partial isometry, then there exists a partial
isometry $\hat{%
v}$ with $\left\| v-\hat{v}\right\| \leq \delta $.
\end{lma}

\begin{proof}
Assume $\epsilon <\frac{1}{4}.$ Lemma 2.3 then provides a projection $%
p\in C^{\ast }(v^{\ast }v)$ with $\left\| v^{\ast }v-p\right\|
<2\epsilon .$
Let $\overline{v}\left| v\right|$ be the polar decomposition for $v$
and note
that each of $p,\;v^{\ast }v,\;\left| v\right| $ and $\overline{v}^{\ast
}\overline{v%
}$ is an element in the abelian $C^{\ast }$-algebra $C^{\ast }(v^{\ast
}v).$ Now let $\hat{v}=\overline{v}p.$ We show that $\hat{v}$
is the required
partial isometry. Firstly, $\hat{v}$ is indeed a
partial isometry: $\hat{v}^{\ast }\hat{v}=(%
\overline{v}p)^{\ast }\overline{v}p=p\overline{v}^{\ast
  }\overline{v}p=pqp$, which, since $p$ and $q$ commute, is itself a
projection. Secondly, we estimate $%
\left\| v-\hat{v}\right\| =\left\| \overline{v}\left| v\right| -%
\overline{v}p\right\| \leq \left\| \overline{v}\right\| \left\| \left|
v\right| -p\right\| \leq \left\| \left| v\right| -p\right\| .$ The proof
is concluded by noting that $%
\left\| v^{\ast }v-p\right\| <2\epsilon $ and the square root map is
continuous on the positive cone of a $C^{*}$-algebra.
\end{proof}
In fact the partial isometry  $\hat{v}$ is dominated by the partial
isometry
$\overline{v}$, that is, $\hat{v}^{\ast }\hat{v}(\mathcal{H})\subseteq
\overline{v}^{\ast }\overline{v}%
(\mathcal{H})$ and $\hat{v}\hat{v}^{\ast }(\mathcal{H})\subseteq
\overline{v}\overline{v}%
^{\ast }(\mathcal{H}).$ To see this, consider:
\begin{equation*}
\overline{v}%
^{\ast }\overline{v}\hat{v}^{\ast }\hat{v}=qpqp=pqp=\hat{v}%
^{\ast }\hat{v}\Rightarrow \hat{v}^{\ast }\hat{v}(\mathcal{H})\subseteq 
\overline{v}^{\ast }\overline{v}(\mathcal{H}) 
\end{equation*}
 \begin{equation*}
\overline{v}\overline{v}^{\ast }\hat{v}\hat{v}^{\ast }=\overline{v}%
\overline{v}^{\ast }\overline{v}p\overline{v}^{\ast
}=\overline{v}\overline{v%
}^{\ast }\overline{v}pp\overline{v}^{\ast
}=\overline{v}p\overline{v}^{\ast }%
\overline{v}p\overline{v}^{\ast }=\left( \hat{v}\hat{v}^{\ast }\right)
^{2}=\hat{v}%
\hat{v}^{\ast }\Rightarrow \hat{v}\hat{v}^{\ast }(\mathcal{H})\subseteq
\overline{v}%
\overline{v}^{\ast }(\mathcal{H}).
\end{equation*}
As a consequence of this, if $v$ and $w$ are orthogonal
approximate partial isometries, then we may infer the existence of
orthogonal  exact partial
isometries $\hat{v},\hat{w}$ close to $v$ and $w$ respectively.
\begin{lma}
Let $v=\left[ 
\begin{array}{cc}
v_{1} & v_{2} \\ 
\varepsilon  & v_{3}
\end{array}
\right] $ be a partial isometry with block diagonal final
projection and suppose $\left\| \varepsilon \right\| \leq \epsilon .$
Then
there exists a partial isometry $\hat{v}=\left[ 
\begin{array}{cc}
\hat{v}_{1} & \hat{v}_{2} \\ 
0 & \hat{v}_{3}
\end{array}
\right] $ with $\left\| v-\hat{v}\right\| \leq \delta (\epsilon
).$
\end{lma}
Note that although the matrix for $v$ is assumed to be square, we make
no
such assumptions for the diagonal sub operators. It is open whether
this result remains true when the block diagonality assumption on the
final projection is dropped. 
\begin{proof}
Firstly, write
$v=\left[ 
\begin{array}{c}
w_{1} \\ 
w_{2}
\end{array}
\right] ,$ where $w_{1}=\left[ 
\begin{array}{cc}
v_{1} & v_{2}
\end{array}
\right] $ and $w_{2}=\left[ 
\begin{array}{cc}
\varepsilon  & v_{3}
\end{array}
\right]$.
\vspace{0.2cm}
By hypothesis, $vv^{\ast }$ is a block diagonal projection:
\begin{equation*}vv^{\ast }=\left[ 
\begin{array}{cc}
v_{1}v_{1}^{\ast }+v_{2}v_{2}^{\ast } & 0 \\ 
0 & \varepsilon \varepsilon ^{\ast }+v_{3}v_{3}^{\ast }
\end{array}
\right] \end{equation*}
from which we see that $ w_{1}w_{1}^{\ast }$ and $w_{2}w_{2}^{\ast }$
are projections. It follows that $w_{1}$ and $w_{2}$ are
partial isometries and since $v$
itself is a partial isometry, $w_{1}$ and $w_{2}$ must be orthogonal.
Now we
claim that $v_{3}$ is an approximate partial isometry. To show this we
estimate $%
\left\| v_{3}v_{3}^{\ast }-\left( v_{3}v_{3}^{\ast }\right) ^{2}\right\|
.$
Since $\varepsilon \varepsilon ^{\ast }+v_{3}v_{3}^{\ast }$ is a
projection, 
$\varepsilon \varepsilon ^{\ast }+v_{3}v_{3}^{\ast }=\left( \varepsilon
\varepsilon ^{\ast }+v_{3}v_{3}^{\ast }\right) ^{2}$ and thus;
\begin{eqnarray*}\left\|
v_{3}v_{3}^{\ast }-\left( v_{3}v_{3}^{\ast }\right) ^{2}\right\| &\leq&
2\left\| v_{3}v_{3}^{\ast }\varepsilon \varepsilon ^{\ast }\right\|
+\left\|
\left( \varepsilon \varepsilon ^{\ast }\right) ^{2}\right\| +\left\|
\varepsilon \varepsilon ^{\ast }\right\| \\&\leq& \epsilon ^{2}(\epsilon
^{2}+3).
\end{eqnarray*}
Provided $\epsilon ^{2}(\epsilon ^{2}+3)<\frac{1}{4}$, Lemma 2.4
ensures the existence of a partial isometry $\hat{v}_{3}$ with $\left\|
v_{3}-\hat{v}%
_{3}\right\| \leq \delta _{1}(\epsilon ).$

Now set $\overline{w}_{2}=\left[ 
\begin{array}{cc}
0 & \hat{v}_{3}
\end{array}
\right] $. We obtain a partial isometry with range orthogonal to $w_1$
and with $\left\| w_{2}-\overline{w}_{2}\right\| \leq \epsilon +\delta
_{1}(\epsilon )$. Set $P=\overline{w}_{2}^{\ast }\overline{w}_{2}$ and
let $%
\widetilde{w}_{1}=w_{1}P^{\perp }.$ We aim to show that
$\widetilde{w}_{1}$
is close to $w_{1}.$ Firstly we claim that $\overline{w}_{2}$ and
$w_{1}$
are $\epsilon +\delta _{1}(\epsilon )$-approximately orthogonal. Since
they
have orthogonal ranges, $\left\| w_{1}^{\ast }\overline{w}_{2}\right\|
=0.$
Now,
\begin{eqnarray*}
\left\| \overline{w}_{2}w_{1}^{\ast }\right\| &=&\left\| (\overline{w}%
_{2}-w_{2})w_{1}^{\ast }+w_{2}w_{1}^{\ast }\right\| \\&\leq& \left\|
\overline{w}%
_{2}-w_{2}\right\| \left\| w_{1}^{\ast }\right\| +\left\|
w_{2}w_{1}^{\ast }\right\| \\&\leq& \epsilon +\delta
_{1}(\epsilon)
\end{eqnarray*}
We now estimate the distance between $\widetilde{w}_{1}$ and $w_{1}$: 
\begin{eqnarray*} 
\left\|\widetilde{w}_{1}-w_{1}\right\| &=&\left\| w_{1}P^{\perp
}-w_{1}\right\|
 \\&=&\left\| w_{1}\overline{w}_{2}^{\ast }\overline{w}_{2}\right\|
\\&\leq& \epsilon
+\delta _{1}(\epsilon ). 
\end{eqnarray*}
Now $\widetilde{w}_{1}$ need not be a partial
isometry but we claim it is close to one. This is immediate  from the
following estimate:
\begin{eqnarray*}
\left\|  \widetilde{w}_{1}\widetilde{w}_{1}^{\ast }-\left(
\widetilde{w}_{1}%
\widetilde{w}_{1}^{\ast }\right) ^{2}\right\| &=&\left\|
\widetilde{w}_{1}%
\widetilde{w}_{1}^{\ast }-w_{1}w_{1}^{\ast }+\left( w_{1}w_{1}^{\ast
}\right) ^{2}-\left( \widetilde{w}_{1}\widetilde{w}_{1}^{\ast }\right)
^{2}\right\| \ \\&\leq& \delta _{2}(\epsilon ),
\end{eqnarray*}
since $w_{1}$ is a partial isometry. Provided $\delta
_{2}(\epsilon )<\frac{1}{4}$, $\widetilde{w}_{1}$ is close to a partial
isometry $\overline{w}_{1}=\left[ 
\begin{array}{cc}
\hat{v}_{1} & \hat{v}_{2}
\end{array}
\right] $ dominated by $\ \widetilde{w}_{1}.$ Thus $\overline{w}_{1}$ is
orthogonal to $\overline{w}_{2}.$ Now let $\hat{v}=\left[ 
\begin{array}{c}
\overline{w}_{1} \\ 
\overline{w}_{2}
\end{array}
\right] $ which by construction is a partial isometry of the required
form.
To complete the proof we estimate:
\begin{eqnarray*}
\left\| v-\hat{v}\right\| &\leq&
\left\| w_{1}-\overline{w}_{1}\right\| +\left\| w_{2}-\overline{w}%
_{2}\right\| \\&\leq& 2\epsilon +2\delta _{1}(\epsilon )+\delta
_{2}(\epsilon).
\end{eqnarray*}
\end{proof}
We shall use this result as the initial step in an induction to
generalise to an $r\times r$ block decomposition. In order to do this
however, we shall need to call on the following elementary lemmas, the
first of
which is a standard $C^{\ast}$-algebra result. For a proof,
see~\cite{murp}.

\begin{lma}
Let $p,q$ be projections in a unital $C^{*}$-algebra $\mathcal{A}$, and
suppose
that $\left\| q-p\right\| <1.$ Then there exists a unitary
$u\in\mathcal{A}$ such
that $q=upu^{\ast }$ and $\left\| I-u\right\| \leq \sqrt{2}\left\|
q-p\right\| ,$ namely, $u=v\left| v\right| ^{-1},$ where $v=I-p-q+2pq.$
\end{lma}
\begin{lma}
Let $w=(w_{i,j})$ be an operator which is $\epsilon _{1}$-close to a
partial isometry 
$v$ and is such that $ww^{\ast }$ is $\epsilon _{2}$-close to a block
diagonal projection. Then $w$ is $\delta (\epsilon_{1},
\epsilon_{2})$-close to a partial
isometry $\overline{w}$ with
$\overline{w}\hspace{0.1cm}\overline{w}^{\ast }$ block
diagonal.
\end{lma}

\begin{proof}
$\left\| w-v\right\| \leq \epsilon _{1}$ from which we infer that
$p^{\prime}=vv^{\ast }$
is an approximately block diagonal projection, that is there exists a
block
diagonal projection $p$ with $\left\| p^{\prime }-p\right\| \leq \delta
(\epsilon_{1},\epsilon_{2}) <1,$ say. Lemma 2.6 now provides a unitary
$u\in M_{n}$
such that $p=up^{\prime }u^{\ast }$ and $\left\| I-u\right\| \leq
\sqrt{2}%
\left\| p^{\prime }-p\right\| \leq \sqrt{2}\delta
(\epsilon_{1},\epsilon_{2}).$ Now set $%
\overline{w}$ =$uv.$ Then $\overline{w}\hspace{0.1cm}\overline{w}^{\ast
}=uvv^{\ast
}u^{\ast }=up^{\prime }u^{\ast }=p,$ a block diagonal projection. Thus
$%
\overline{w}$ $\ $is a partial isometry with the required properties
and
\begin{eqnarray*} 
\left\| w-\overline{w}\right\| &\leq& \left\| w-v\right\| +\left\|
v-\overline{w} \right\| \\&\leq& \epsilon _{1}+\left\| v-uv\right\|
\\&\leq& \epsilon
_{1}+\left\| I-u\right\| \left\| v\right\| \\&\leq& \epsilon
_{1}+\sqrt{2}
\delta (\epsilon_{1},\epsilon_{2}).
\end{eqnarray*}
\end{proof}

\begin{lma}
Let $v,w\in M_{n}$ be partial isometries with $\left\| v-w\right\| <1.$
Then 
$v$ and $w$ have equal rank.
\end{lma}
See~\cite{hal} for a proof. We now have all the necessary machinery in
place to prove the most
important result on our way to the proof of Theorem 2.1, which in
essence  is a non-selfadjoint generalisation of Lemma 2.4. 
\begin{lma}
Let $v=\left[ 
\begin{array}{cccc}
v_{1,1} & v_{1,2} & \cdots & v_{1,r} \\ 
\varepsilon _{2,1} & v_{2,2} & \cdots & v_{2,r} \\ 
\vdots & \vdots & \ddots & \vdots \\ 
\varepsilon _{r,1} & \varepsilon _{r,2} & \cdots & v_{r,r}
\end{array}
\right] $ be a partial isometry with block diagonal final projection,
in which the operators below the main block diagonal
have norm $\leq \epsilon .$ Then there exists a positive function
$\delta
(\epsilon ,r),$ with $\delta(\epsilon ,r)\rightarrow0$ as
$\epsilon\rightarrow0$ for fixed $r$ 
(and which is independent of the sizes of the
blocks of $v$) such that for each $\epsilon$ there exists a block
upper triangular partial isometry
$\hat{v}=\left[ 
\begin{array}{cccc}
\hat{v}_{1,1} & \hat{v}_{1,2} & \cdots & \hat{v}_{1,r} \\ 
0 & \hat{v}_{2,2} & \cdots & \hat{v}_{2,r} \\ 
\vdots & \vdots & \ddots & \vdots \\ 
0 & 0 & \cdots & \hat{v}_{r,r}
\end{array}
\right] $
with $\left\| v-\hat{v}\right\| \leq \delta
(\epsilon ,r).$
\end{lma}

\begin{proof}
 Lemma 2.5 gives the result for $r=2$, so we assume the result to
hold for all $v\in M_{n}$ satisfying the conditions with $s\times s$
block
decomposition for $s<r$ and consider $v$ as in the statement. Write
$v=\left[ 
\begin{array}{cc}
v_{1,1} & X \\ 
\varepsilon  & U
\end{array}
\right] $
with the obvious identifications. We shall also refer to the
decomposition
$v=\left[ 
\begin{array}{c}
w_{1} \\ 
w_{2}
\end{array}
\right] ,$ where, $w_{1}=\left[ 
\begin{array}{cc}
v_{1,1} & X
\end{array}
\right] $ and $w_{2}=\left[ 
\begin{array}{cc}
\varepsilon  & U
\end{array}
\right] .$
As in the proof of Lemma 2.5, $vv^{\ast }$ is block diagonal, implying
that $%
w_{1}$and $w_{2}$ are orthogonal partial isometries, a fact we make
use of later. Now,
\begin{equation*}
vv^{\ast }=\left[ 
\begin{array}{cc}
v_{1,1}v_{1,1}^{\ast }+XX^{\ast } & v_{1,1}\varepsilon ^{\ast }+XU^{\ast
}
\\ 
\left( v_{1,1}\varepsilon ^{\ast }+XU^{\ast }\right) ^{\ast } &
\varepsilon
\varepsilon ^{\ast }+UU^{\ast }
\end{array}
\right] 
\end{equation*}
and the block diagonality assumption yields
\begin{equation*}
v_{1,1}\varepsilon ^{\ast }+XU^{\ast }=\left( v_{1,1}\varepsilon ^{\ast
}+XU^{\ast }\right) ^{\ast }=0.
\end{equation*}
Thus $\varepsilon \varepsilon ^{\ast }+UU^{\ast }$ is a block diagonal
projection. Since $\left\| \varepsilon \right\| $ is small, we infer
that $U$ is  close to a partial isometry. Also,
clearly $UU^{\ast }$ is approximately block diagonal, from which  Lemma
2.7
guarantees the existence of a partial isometry $U^{\prime }$ with
block diagonal final projection  and with $\left\| U-U^{\prime }\right\|
$ $\leq
\delta _{1}(\epsilon ).$ Now
\begin{equation*}
U^{\prime }=\left[ 
\begin{array}{ccc}
U_{2,2}^{\prime } & \cdots & U_{2,r}^{\prime } \\ 
\vdots & \ddots & \vdots \\ 
U_{r,2}^{\prime } & \cdots & U_{r,r}^{\prime }
\end{array}
\right] ,
\end{equation*} 
with the operators below the block diagonal having norm $\leq \delta
_{1}(\epsilon )+\epsilon$. We now apply the induction hypothesis to
$U^{\prime }$ to
discover a block upper triangular partial isometry $\hat{U}$ close to
$U^{\prime }$, within $%
\delta _{2}(\epsilon )$ say.
Let
\begin{equation*}
v^{\prime }=\left[ 
\begin{array}{cc}
v_{1,1} & X \\ 
0 & \hat{U}
\end{array}
\right] =\left[ 
\begin{array}{c}
w_{1} \\ 
w_{2}^{\prime }
\end{array}
\right] .
\end{equation*}
Now $v^{\prime }$ need not be a partial isometry, since $w_{1}$ and $%
w_{2}^{\prime }$ need not be orthogonal. However they are approximately
orthogonal:
\begin{eqnarray*} \left\| w_{2}^{\prime }w^{\ast}_{1}\right\| &=&\left\|
(w_{2}^{\prime
}-w_{2})w_{1}^{\ast }+w_{2}w_{1}^{\ast }\right\| \\&\leq& (r-1)\epsilon
+\delta
_{1}(\epsilon )+\delta _{2}(\epsilon ),  
\end{eqnarray*}
since $\left\| \varepsilon \right\| \leq (r-1)\epsilon $, 
$\left\| \hat{U}-U\right\| \leq \delta _{1}(\epsilon )+\delta
_{2}(\epsilon )$ and $w_{1}\bot w_{2}.$ 
Let $\widetilde{w}_{1}=w_{1}Q^{\bot }$ where $Q=\left(
w_{2}^{\prime }\right) ^{\ast }w_{2}^{\prime }.$ Then
$\widetilde{w}_{1}$and 
$w_{2}^{\prime }$ are orthogonal. Now $\widetilde{w}_{1}$ need not be a
partial isometry, but it is close to one which it dominates, 
$w_{1}^{\prime }$,  with
$%
\left\| \widetilde{w}_{1}-w_{1}^{\prime }\right\| \leq \delta
_{3}(\epsilon )
$ say.
This follows since $\widetilde{w}_{1}$ is clearly an approximate partial
isometry.  Now
\begin{equation*}
\hat{v}=\left[ 
\begin{array}{c}
w_{1}^{\prime } \\ 
w_{2}^{\prime }
\end{array}
\right]
\end{equation*}
is a partial isometry by construction with the correct form and

\begin{eqnarray*}
\left\| v-\hat{v}\right\| &\leq& \left\| w_{1}-w_{1}^{\prime }\right\|
+\left\| w_{2}-w_{2}^{\prime }\right\| \\
&\leq& \left\| w_{1}-\widetilde{w}_{1}\right\| +\left\| \widetilde{w}
_{1}-w_{1}^{\prime }\right\| +\left\| w_{2}-w_{2}^{\prime } \right\| \\
&\leq& 2(r-1)\epsilon +2\delta _{1}(\epsilon )+2\delta _{2}(\epsilon
)+\delta
_{3}(\epsilon)
\end{eqnarray*}
The result follows by induction.
\end{proof}
We remark that, as a consequence of Lemma 2.8, if $v$ is unitary then
$\hat{v}$ must also be unitary
provided $\left\| v-\hat{v}\right\| $ is small enough.
\begin{lma}
Let $I_{N}=I_{n_{1}}\oplus I_{n_{2}}\oplus ...\oplus I_{n_{r}}$ be a
decomposition of the identity operator in $M_{N}$, imposing a block
structure. Take $P,Q$ to be orthogonal block diagonal projections of the
same
rank and let $b\in M_{N}$ be a partial isometry such that $b^{\ast
}b=Q,\;bb^{\ast }=P$ and $\left\| I_{n_{i}}bI_{n_{j}}\right\| \leq
\epsilon $
for $i>j.$ Then there exists a partial isometry $\hat{b}\in M_{N}$ with:
\begin{enumerate}
\item$I_{n_{i}}\hat{b}I_{n_{j}}=0\;$for $i>j,$
\item$\hat{b}^{\ast }\hat{b}=Q\;$and $\hat{b}\hspace{0.1cm}\hat{b}^{\ast
}=P,$
\item$\left\| b-\hat{b}\right\| \leq \delta (\epsilon ,r).$
\end{enumerate}
\end{lma}
\begin{proof}
Let $P_{i}=I_{n_{i}}P,\;Q_{i}=I_{n_{i}}Q$ and let
$\mathcal{H}_{i}=$ran$P_{i},$ $%
\mathcal{K}_{i}=$ran$Q_{i}$ for $1\leq i\leq r.$ Set
$\mathcal{H}=\mathcal{H}_{1}\oplus ...\oplus
\mathcal{H}_{r},\;\mathcal{K}=\mathcal{K}_{1}\oplus ...\oplus
\mathcal{K}_{r}$. Then $\dim \mathcal{H}=\dim \mathcal{K}$ and $%
b:\mathcal{K}\rightarrow \mathcal{H}$ is a unitary operator with reduced
matrix
\begin{equation*}
\left[ 
\begin{array}{cccc}
b_{1,1} & b_{1,2} & \cdots  & b_{1,r} \\ 
b_{2,1} & b_{2,2} & \cdots & b_{2,r} \\ 
\vdots & \vdots & \ddots & \vdots \\ 
b_{r,1} & b_{r,2} & \cdots & b_{r,r}
\end{array}
\right] \in M_{S}
\end{equation*}
 where $s=rk(P)=rk(Q).$ Now
$b_{i,j}:\mathcal{K}_{j}\rightarrow \mathcal{H}_{i}$ and
$b_{i,j}=P_{i}bQ_{j}$ with $\left\|
b_{i,j}\right\| \leq \epsilon \;$for all$\;i>j.$ Lemma 2.9  now supplies
a unitary operator $\hat{b}$ with matrix
\begin{equation*}
\left[ 
\begin{array}{cccc}
\hat{b}_{1,1} & \hat{b}_{1,2} & \cdots & \hat{b}_{1,r} \\ 
0 & \hat{b}_{2,2} & \cdots & \hat{b}_{2,r} \\ 
\vdots & \vdots & \ddots & \vdots \\ 
0 & 0 & \cdots & \hat{b}_{r,r}
\end{array}
\right]  
\end{equation*}
This operator is upper triangular with respect to the same decomposition
and with $\left\| b-\hat{b}%
\right\| \leq \delta (\epsilon ,r).$ Since $\hat{b}$ is unitary,
$({\hat{b}})^{\ast }\hat{b%
}=Q\;$and $\hat{b}\hspace{0.1cm}(\hat{b})^{\ast }=P.$ Now in the
original matrix for 
$b$ replace each $I_{n_{i}}bI_{n_{j}}$ by $\hat{b}_{i,j}$ for
$i\leq j$ and  $0$ for $i>j$ to find the required partial isometry.
\end{proof}
The proof of our next Lemma uses standard techniques and is an obvious
extension of the proof of the similar Lemma 3.5 in~\cite{hef} and as
such
we shall omit the details. 
\begin{lma}
Let $\mathcal{A}_{1}$,$\mathcal{A}_{2}\in\mathcal{N}_r$.  Suppose that
$\epsilon >0$ and $\mathcal{A}_1%
\subseteq_{\epsilon}\mathcal{A}_2$. Then $(\mathcal{A}_{1}\cap
\mathcal{A}_{1}^{\ast })
\subseteq_{\delta (\epsilon ,r)}(\mathcal{A}_{2}\cap
\mathcal{A}_{2}^{\ast })$.
\end{lma}
We end this section with the proof of Theorem 1.
\begin{proof}
Let $\mathcal{A}_{1}$ and $\mathcal{A}_{2}$ be nest algebras each of
which is determined
by $r$ non-zero nest projections and suppose  $\mathcal{A}_1%
\subseteq _{\epsilon }\mathcal{A}_2$.
Denote the matrix unit
system for $\mathcal{A}_{1}$ by $\left\{ e_{i,j}\right\} _{i,j\in
I}$. We wish to create a matrix unit
system, $ \left\{ f_{i,j}\right\} _{i,j\in I}$   for a copy of
$\mathcal{A}_{1}$ inside  $\mathcal{A}_{2},$
satisfying $\left\| e_{i,j}-f_{i,j}\right\| \leq \delta (\epsilon).$ We
then set the star extendible
map $\phi $ to be the linear extension of the correspondences $%
e_{i,j}\rightarrow f_{i,j}$ for each $i,j.$ By hypothesis and Lemma
2.11, $%
\mathcal{A}_{1}\cap \mathcal{A}_{1}^{\ast }\subseteq _{\delta
_{1}(\epsilon
  )} %
\mathcal{A}_{2}\cap \mathcal{A}_{2}^{\ast },$ and it is a standard self
adjoint result that this
ensures the existence of a matrix unit system for a copy of
$\mathcal{A}_{1}\cap
\mathcal{A}_{1}^{\ast }$ in $\mathcal{A}_{2}\cap \mathcal{A}_{2}^{\ast
},$ $\left\{ f_{i,j}\right\}
_{i,j\in J}$, satisfying $\left\| e_{i,j}-f_{i,j}\right\| \leq
\delta _{2}(\epsilon )$ for all$\;i,j\in J.$

Now since $\mathcal{A}_{1}$ is a nest algebra, its reduced digraph is
generated by a
bilateral tree and every
non-selfadjoint matrix unit may be factorized uniquely into a product of
matrix units associated with the tree. Label the matrix units
corresponding to the edges of the tree by  $e_{1},...,e_{r-1}.$ We need
to
find partial isometries $f_{1},...,f_{r-1}$ close to each
$e_{1},...,e_{r-1}$
respectively and having the correct initial and final projections, that
is,
if $e_{1}^{\ast }e_{1}=e_{p,p}$ and $e_{1}e_{1}^{\ast }=e_{q,q}$, then
we
demand that $f_{1}^{\ast }f_{1}=f_{p,p}$ and $f_{1}f_{1}^{\ast
}=f_{q,q}.$
Having found suitable partial isometries as described, we may define the
remaining matrix units without ambiguity by multiplication. More
precisely
suppose that the matrix unit $e_{1}$ satisfies
$e_{1}=I_{n_{k}}e_{1}I_{n_{l}}$
(so that $e_{1}$ lives in the $k,l$ block) and suppose that we have
found a
partial isometry $f_{1}\in \mathcal{A}_{2}$ with the correct initial and
final
projections as described above and with $\left\| e_{1}-f_{1}\right\|
\leq
\delta _{3}(\epsilon ).$ Then for
$n_{1}+\cdots+n_{k-1}+1\leq i\leq n_{k}$, and
$n_{1}+\cdots+n_{l-1}+1\leq j\leq
n_{l}$
(that is $i,j$ corresponding to the $k,l$ block) we define $%
f_{i,j}=f_{i,q}f_{1}f_{p,j}.$ Then for such $i,j$ we have 
\begin{eqnarray*}
\left\|
e_{i,j}-f_{i,j}\right\| &\leq& \left\| e_{i,q}-f_{i,q}\right\| +\left\|
e_{1}-f_{1}\right\| +\left\| e_{p,j}-f_{p,j}\right\| \\&\leq& 2\delta
_{2}(\epsilon )+\delta _{3}(\epsilon ). 
\end{eqnarray*}
Having generated all the matrix units
we can using the matrix units $\left\{ f_{i,j}\right\} _{i,j\in J}$ and
the
partial isometries $f_{1},...,f_{r-1}$, the remaining matrix units are
determined in a well defined way. Thus a full matrix unit system for a
copy
of $\mathcal{A}_{1}$ in $\mathcal{A}_{2}$ has been created. Note that
the quantity $\underset{%
i,j\in I}{\sup }\left\| e_{i,j}-f_{i,j}\right\| $ depends only on
$\epsilon $
and $r$; the error does not increase within a single block. The
remainder of the proof shall be concerned with finding suitable
candidates
for $f_{1},...,f_{r-1}.$ We demonstrate the technique only for $f_{1}$
as
the other cases are similar.

Consider $e_{1}\in ball(\mathcal{A}_{1})$. We keep the same notation as
above for its
projections$.$ By hypothesis there exists $b_{1}\in
ball(\mathcal{A}_{2})$ with $%
\left\| e_{1}-b_{1}\right\| \leq \epsilon .$ Let $P=f_{p,p}$ ,
$Q=f_{q,q}.$ We
use these projections to cut down $b_{1}.$
Set $b=Qb_{1}P$ so that $b=QbP,$ and we have  $b^{\ast }b\in PM_{n}P$ \
(where $n=n_{1}+\cdots+n_{r}$). By the polar decomposition, write
$b=\hat{b}%
\left| b\right| .$ Then 
\begin{eqnarray*}
\left\| e_{1}-b\right\| &=&\left\| e_{1}-QbP\right\| \\
&\leq& 3{\delta _{2}}(\epsilon ), 
\end{eqnarray*}
since, $\left\| e_{1}^{\ast }e_{1}-P\right\| \leq \delta _{2}(\epsilon
),
$ $\left\| e_{1}e_{1}^{\ast }-Q\right\| \leq \delta _{2}(\epsilon )$ and
$%
\left\| e_{1}-b_{1}\right\| \leq \epsilon .$
Also we have 
\begin{eqnarray*}
\left\| b^{\ast }b-P\right\| &\leq& \left\| b^{\ast }b-e_{1}^{\ast
}e_{1}\right\| +\left\| e_{1}^{\ast }e_{1}-P\right\| 
 \\&\leq& \left\| (b^{\ast }-e_{1}^{\ast })b+e_{1}(b-e_{1})\right\|
+\left\|
e_{1}^{\ast }e_{1}-P\right\| 
 \\&\leq& 7\delta_{2}(\epsilon)
\end{eqnarray*}
Provided $\epsilon $ was originally small enough so that     $
7\delta_{2}(\epsilon)<$ 1, we have that $b^{\ast }b$ is invertible in $%
PM_{n}P$, from which  $(b^{\ast }b)^{\frac{1}{2}}$ is invertible in
$PM_{n}P$.
Thus $(\hat{b})^{\ast }\hat{b}=P.$ Similarly,
$\hat{b}\hspace{0.1cm}(\hat{b})%
^{\ast }=Q.$ We now show that $b$ and $\hat{b}$ are close and hence that

$e_{1}$ and $\hat{b}$ are close: 
\begin{eqnarray*}
\left\| b-\hat{b}\right\| &=&\left\| 
\hat{b}\left| b\right| -\hat{b}\right\| \\&=&\left\| \hat{b}(\left|
b\right| -P)\right\|  \\&\leq& \left\| \left| b\right| -P\right\|
\\&\leq& \left\|
b^{\ast }b-P\right\| \\&\leq&  7\delta_{2}(\epsilon)
\end{eqnarray*}
We see that $\hat{b}$ is a partial isometry with $(\hat{b})^{\ast }%
\hat{b}=P$ and $\hat{b}\hspace{0.1cm}(\hat{b})^{\ast}=Q$ (which are
block
diagonal) and so would be a candidate for our required partial isometry
$%
f_{1}$, but we cannot guarantee that $f_{1}\in \mathcal{A}_{2}.$ However
$f_{1}$ is a
partial isometry in $M_{n}$ and since it is within $  
7\delta_{2}(\epsilon)    $ of $b$ which is in $\mathcal{A}_{2}$, it must
be the case that the
operators in the zero blocks for $\mathcal{A}_{2}$ in the representing
matrix of $f_{1}
$ have norm $\leq   7\delta_{2}(\epsilon).$ Lemma 2.10 now
provides a partial isometry $v\in \mathcal{A}_{2}$ with the same initial
and final projections as 
$\hat{b}$ and with $\left\| \hat{b}-v\right\| \leq \delta
_{3}(\epsilon).$
Setting $f_{1}=v$ completes the proof.
\end{proof}

\section{Characterisation of regular limit algebras}

In this section, we continue the theme of the local characterisation of
separable operator algebras, but in the regular setting. The imposition
of
regularity allows us greater freedom to
characterise limits from a much wider class of algebra than
previously, in terms of a local description capturing this additional
structure. As we have seen, the difficulty in characterising general
star extendible limits formed from arbitrary non-selfadjoint families
lies in showing that such a family is stable. In the regular context,
we need to overcome the analogous problem of showing that a family
whose regular limits we wish to characterise, is regularly stable. In
essence, this means that an approximate regular containment of one
family member in another is close to an exact regular containment.
  
Throughout, we shall assume our operator algebras to
contain a masa of the generated $C^{\ast }$-algebra. We shall denote the
pair of operator algebra $\mathcal{A}$ and masa $\mathcal{C}$ by
$(\mathcal{A},\mathcal{C}),$ and take a chain $%
(\mathcal{A}_{1},\mathcal{C}_{1})\subseteq
(\mathcal{A}_{2},\mathcal{C}_{2})\subseteq \cdots$ to mean
$\mathcal{A}_{k}\subseteq
\mathcal{A}_{k+1},\mathcal{C}_{k}\subseteq \mathcal{C}_{k+1}$ for each
$k.$
Furthermore, the set of $\mathcal{C}$ normalising partial isometries
in $\mathcal{A}$, that is, partial isometries $v$ having the property
that $c \in \mathcal{C}$ implies that $v^{*}cv$ and $vcv^{*}$ both
belong to $\mathcal{C}$, will be denoted by the symbol
$N_{\mathcal{C}}(\mathcal{A}).$ 
Whenever we speak of a
regular injection, $\phi ,$ between algebras
$(\mathcal{A}_{1},\mathcal{C}_{1})$ and $%
(\mathcal{A}_{2},\mathcal{C}_{2}),$ it will be taken to mean that the
map is regular with respect to the
masas written, that is $\phi :\mathcal{A}_{1}\rightarrow
\mathcal{A}_{2}$ and $\phi
(N_{\mathcal{C}_{1}}(\mathcal{A}_{1}))\subseteq
N_{\mathcal{C}_{2}}(\mathcal{A}_{2}).$

\begin{defn}
Let $\mathcal{A}$ be an operator algebra with masa $\mathcal{C}.$ An
element $w\in \mathcal{A}$ is said to
be $\epsilon $-approximately normalising if for each $c\in \mathcal{C}$
with $\left\|
c\right\| \leq 1$ there exist $c_{1},$ $c_{2}\in \mathcal{C}$ with
$\left\| w^{\ast
}cw-c_{1}\right\| <\epsilon $ and $\left\| wcw^{\ast }-c_{2}\right\|
<\epsilon .$
\end{defn}

\begin{defn}
Let $(\mathcal{A}_{1},\mathcal{C}_{1})$ and
$(\mathcal{A}_{2},\mathcal{C}_{2})$ be operator algebras with masas. We
shall write $\ N_{\mathcal{C}_{1}}(\mathcal{A}_{1})\subseteq_{\epsilon}%
N_{\mathcal{C}_{2}}(\mathcal{A}_{2})$ to mean for every $v_{1}\in
N_{\mathcal{C}_{1}}(\mathcal{A}_{1})$   
there exists $%
v_{2}\in N_{\mathcal{C}_{2}}(\mathcal{A}_{2})$ with $\left\|
v_{1}-v_{2}\right\| <\epsilon .$
\end{defn}

\begin{lma}
Let $w\in ball(\mathcal{A})$ be an $\epsilon $-approximately
normalising, $\epsilon $-approximate
partial isometry. Then for each projection $p\in proj(\mathcal{C}),$
there exist
projections $p_{1},$ $p_{2}\in proj(\mathcal{C})$ such that $\left\|
w^{\ast
}pw-p_{1}\right\| ,$ $\left\| wpw^{\ast }-p_{2}\right\| <\delta
(\epsilon
)$
\end{lma}

\begin{proof}
Firstly we show that for any $p\in proj(\mathcal{C}),\;wpw^{\ast }$ is
an approximate
projection.  By hypothesis, there exists $c\in \mathcal{C}$ with 
$\left\| w^{\ast }w-c\right\| <\epsilon .$ Then 
\begin{eqnarray*}
\left\| w^{\ast
}wp-pw^{\ast }w\right\| &=&\left\| w^{\ast }wp-cp+pc-pw^{\ast
}w\right\|  \\&\leq&
\left\| w^{\ast }w-c\right\| \left\| p\right\| +\left\| p\right\|
\left\|
c-w^{\ast }w\right\|  \\&<&2\epsilon .
\end{eqnarray*}
Thus $w^{\ast }w$ approximately commutes with any $p\in
proj(\mathcal{C}).$ Next,
since $w$ is an $\epsilon$-approximate partial isometry, we show
$ww^{\ast }w\cong w.$
\begin{equation*}
\left\| w\left( w^{\ast }w\right) ^{2}-ww^{\ast }w\right\| \leq \left\|
w\right\| \left\| \left( w^{\ast }w\right) ^{2}-w^{\ast }w\right\|
<\epsilon
.\end{equation*}
Then we have,
\begin{equation*}
\left\| ww^{\ast }w-w\right\| ^{2}=\left\| (ww^{\ast }w-w)^{\ast }\left(
ww^{\ast }w-w\right) \right\| \leq 2\epsilon .
\end{equation*}
Thus, $\left\| ww^{\ast }w-w\right\| \leq (2\epsilon )^{\frac{1}{2}}.$
We are now ready to estimate
\begin{eqnarray*}
\left\| \left( wpw^{\ast }\right) ^{2}-wpw^{\ast }\right\| &=&\left\|
wpw^{\ast }wpw^{\ast }-ww^{\ast }wpw^{\ast }+ww^{\ast }wpw^{\ast
}-wpw^{\ast
}\right\| \\
&\leq& \left\| w\right\| \left\| pw^{\ast }wpw^{\ast }-w^{\ast
}wp\right\|
\left\| w\right\| +\left\| ww^{\ast }w-w\right\| \left\| p\right\|
\left\|
w\right\| \\&<&2\epsilon +(2\epsilon )^{\frac{1}{2}}.
\end{eqnarray*}
Thus for each $p\in proj(\mathcal{C}),$ $wpw^{\ast }$ is an
approximate projection and 
similarly for $w^{\ast }pw.$ Now, $wpw^{\ast }$ is $\epsilon $-close to
an
element $c$ of $\mathcal{C},$ which clearly is an approximate projection
in a $C^{\ast
}$-algebra and so is close to a projection $p_{1}\in \mathcal{C},$
 with $\left\|
p_{1}-c\right\| <\delta _{1}(\epsilon ).$ Then clearly we have
\begin{equation*}
\left\| wpw^{\ast }-p_{1}\right\| \leq \epsilon +\delta _{1}(\epsilon)
\end{equation*}
and the proof is complete.
\end{proof}

The following proposition is a finite dimensional version of the more
general Theorem 9.6 ~\cite{kd} and follows from Arveson's distance
formula.

\begin{prop}
Let $\mathcal{A}$ be a digraph algebra with masa $\mathcal{C}.$ Suppose
$%
w\in \mathcal{A}$ is such that $\left\| wp-pw\right\| <\epsilon $ for
each $p\in
proj(\mathcal{C}).$ Then $dist(w,\mathcal{C})<\epsilon .$
\end{prop}
See~\cite{kd} for a proof.
\begin{prop}
Let $(\mathcal{A},\mathcal{C})$ be a digraph algebra. Let $v\in
ball(\mathcal{A})$ be an $\epsilon$-approximately  normalising
$\epsilon$-approximate partial isometry. Then there exists $
v\in N_{\mathcal{C}}(\mathcal{A})$ with $\left\| \hat{v}-v\right\|
<\delta (\epsilon
).$
\end{prop}

\begin{proof}
The approximately normalising hypothesis easily yields the following
facts. Each entry of $v$ is close in modulus to 1 or 0. If $%
\left| v_{k,l}\right| $ is close to 1, then every other entry in row $k$
and
column $l$ is close to 0. If row $i$ is such that $v_{i}v_{i}^{\ast
}\approx
1$ then there exists a unique index $j$ for which $\left| v_{i,j}\right|
$
is close to 1.
We now define $\mathcal{S}$ to be the index set:
\begin{equation*}
\mathcal{S}=\left\{ (i,j):\left| e_{i,i}ve_{j,j}\right|
  >\frac{1}{2}\right\}
\end{equation*}
and define $v_{1}=\underset{(i,j)\in S}{\sum }e_{i,j}.$
The summary above implies that $v_{1}$ is a partial isometry. 
By Lemma 3.3, for each projection $p\in \mathcal{C}$ there exist
$p^{\prime
},p^{\prime \prime }\in proj(\mathcal{C})$ with $\left\| vpv^{\ast
}-p^{\prime
}\right\| $ and $\left\| v^{\ast }pv-p^{\prime \prime }\right\| $ both
less
than $\delta _{1}(\epsilon ).$ We
assume that $\epsilon $ is sufficiently small so that $\delta
_{1}(\epsilon
)<\frac{1}{4}.$ If $p$ is a minimal projection, then since $vpv^{\ast }$
has
rank 1 or 0 and since $p^{\prime }$ is a projection in $\mathcal{C}$,
$p^{\prime }$
must either be minimal or zero, else $\left\| vpv^{\ast }-p^{\prime
}\right\| \geq 1.$ 

Also, we note
that if $\left| e_{i,i}ve_{j,j}\right| >\frac{1}{\sqrt{2}},$ then
$\left|
e_{j,j}v^{\ast }e_{i,i}ve_{j,j}\right| >\frac{1}{2}$ which implies
$e_{j,j}$
is the minimal projection $p^{\prime }$ satisfying $\left\| v^{\ast
}e_{i,i}v-p^{\prime }\right\| <\epsilon .$ Thus, if $\left\| vpv^{\ast
}-p^{\prime }\right\| <\epsilon ,$ we have, by choice of $v_{1},$ $%
v_{1}^{\ast }p^{\prime }=pv_{1}^{\ast }.$

We now make use of the estimates $\left\| vv^{\ast }v-v\right\| <%
\sqrt{2\epsilon }$ and $\left\| pv^{\ast }v-v^{\ast }vp\right\|
<2\epsilon $
to estimate

\begin{eqnarray*}
\left\| vp-p^{\prime }v\right\|&\leq& \left\| vp-vpv^{\ast }v\right\|
+\left\| vpv^{\ast }v-p^{\prime }v\right\| \\
&\leq& \left\| vp-vv^{\ast }vp\right\| +\left\| vv^{\ast }vp-vpv^{\ast
}v\right\| +\left\| vpv^{\ast }v-p^{\prime }v\right\| \\
&\leq& \left\| v-vv^{\ast }v\right\| \left\| p\right\| +\left\|
  v\right\|
\left\| v^{\ast }vp-pv^{\ast }v\right\| +\left\| vpv^{\ast }-p^{\prime
}\right\| \left\| v\right\| \\
&<&\sqrt{2\epsilon }+3\epsilon \\
&=&\delta _{2}(\epsilon ). 
\end{eqnarray*}
Consider $v_{1}^{\ast }v.$ We have
\begin{eqnarray*}
\left\| (1-p)v_{1}^{\ast }vp\right\| &=&\left\| (1-p)v_{1}^{\ast
}vp-(1-p)v_{1}^{\ast }p^{\prime }v+(1-p)v_{1}^{\ast }p^{\prime
}v\right\| \\
&\leq& \left\| 1-p\right\| \left\| v_{1}^{\ast }\right\| \left\|
vp-p^{\prime
}v\right\| +\left\| (1-p)pv_{1}^{\ast }v\right\| \\ &<&\delta
_{2}(\epsilon),
\end{eqnarray*}
and similarly for $\left\| pv_{1}^{\ast }v(1-p)\right\| .$
Now
\begin{equation*}
\left\| pv_{1}^{\ast }v-v_{1}^{\ast }vp\right\| =\left\| pv_{1}^{\ast
}v(1-p)-(1-p)v_{1}^{\ast }vp\right\| <2\delta _{2}(\epsilon ).
\end{equation*}
Lemma 3.4 now
provides $d\in \mathcal{C}$ with $\left\| v_{1}^{\ast }v-d\right\|
<2\delta
_{2}(\epsilon ).$ Next we need to estimate $\left\| v_{1}v_{1}^{\ast
}v-v\right\| .$ Now $v_{1}v_{1}^{\ast }v$ is $v$ with the rows $v_{i}$
for
which $\left| v_{i}v_{i}^{\ast }\right| <\delta _{1}(\epsilon )$
removed.
Since $vv^{\ast }$ is $\delta _{1}(\epsilon )$ close to a 0,1
projection, we
have $\left\| v_{1}v_{1}^{\ast }v-v\right\| <\delta _{1}(\epsilon
)^{\frac{1%
}{2}}.$ We use this fact firstly to show that $v_{1}^{\ast }v$, and
therefore 
$d$, is an approximate partial isometry. Since $\left\| v^{\ast
}v_{1}v_{1}^{\ast
}v-v^{\ast }v\right\| <\delta _{1}(\epsilon )^{\frac{1}{2}}$  $%
(v_{1}^{\ast }v)^{\ast }v_{1}^{\ast }v$ is a $\delta _{3}(\epsilon
)$-approximate projection. Thus $d$ is an approximate partial isometry
in the
abelian $C^{\ast }$-algebra $\mathcal{C}.$ Replace $d$ by a nearby
partial isometry $%
d^{\prime }\in C$ with $\left\| d-d^{\prime }\right\| <\delta
_{4}(\epsilon)$, so that
\begin{equation*}
\left\| v_{1}^{\ast }v-d^{\prime }\right\| \leq \left\| v_{1}^{\ast
}v-d\right\| +\left\| d-d^{\prime }\right\| <2\delta _{2}(\epsilon)
+\delta _{4 }(\epsilon).
\end{equation*}
Then we have:
\begin{eqnarray*}
\left\| v-v_{1}d\right\| &=&\left\| v-v_{1}v_{1}^{\ast
}v+v_{1}v_{1}^{\ast
}v-v_{1}d\right\| \\ 
&\leq& \delta _{1}(\epsilon)
^{\frac{1}{2}}+2\delta _{2}(\epsilon)
+\delta
_{4}(\epsilon)
=\delta (\epsilon ).
\end{eqnarray*}
Noting that $v_{1}d\in N_{\mathcal{C}}(\mathcal{A})$ concludes the
proof.
\end{proof}
We now introduce some terminology which will be used repeatedly in the
following proofs. Let $(\mathcal{A},\mathcal{C})$ be a digraph algebra
and let $\mathcal{B}$ be the
containing $M_{n}.$ Denote the matrix unit system for $\mathcal{A}$
compatible with
the masa by $\left\{ e_{i,j}\right\} _{i,j}.$ The expectation map $%
E:\mathcal{B}\rightarrow \mathcal{C}$ is defined by
$E(b)=\underset{i}{\sum }e_{i,i}be_{i,i}.$
We note that $\left\| E\right\| =1$ and that
$E(\mathcal{C})=\mathcal{C}.$
\begin{lma}
Let ($\mathcal{A},\mathcal{C})$ be the usual pair of operator algebra
and masa. Let $%
(\mathcal{A}_{1},\mathcal{C}_{1}),$ $(\mathcal{A}_{2},\mathcal{C}_{2})$
be digraph algebras and
masas such that $N_{\mathcal{C}_{i}}(\mathcal{A}_{i})\subseteq
N_{\mathcal{C}}(\mathcal{A})$, $i=1,2$ and $\mathcal{A}_{1}%
\subseteq_{\epsilon}\mathcal{A}_{2},$ where $\epsilon <\frac{1}{4}.$
Then $%
\mathcal{C}_{1}\subseteq \mathcal{C}_{2}.$
\end{lma}
\begin{proof}
Suppose to the contrary that $\mathcal{C}_{1}\nsubseteqq
\mathcal{C}_{2}.$ Since
each $\mathcal{C}_{i}\subseteq \mathcal{C},$ $C^{\ast
}(\mathcal{C}_{1},\mathcal{C}_{2})$ is an abelian $C^{*}$-algebra, and
by the functional calculus, we can find a projection $p\in
\mathcal{C}_{1}\setminus \mathcal{C}_{2},$ with (necessarily),
dist$(p,\mathcal{C}_{2})=1.$ Since $\mathcal{A}_{1}%
\subseteq_{\epsilon}\mathcal{A}_{2}$ there exists $a_{2}\in
ball(\mathcal{A}_{2})$
with $\left\| p-a_{2}\right\| <\epsilon .$ Now consider the expectation
map $%
E_{2}:\mathcal{A}_{2}\rightarrow \mathcal{C}_{2},$
$E_{2}(p)=\underset{j}{\sum }%
e_{j}^{2}pe_{j}^{2}=\underset{j}{\sum }e_{j}^{2}p=I_{2}p,$ where the $%
e_{j}^{2}$ 's are the minimal projections in $\mathcal{C}_{2}$ and
$I_{2}=$ $\underset{%
j}{\sum }e_{j}^{2}$ is the identity in $\mathcal{A}_{2}.$ Then we have
\begin{eqnarray*}
\left\| E_{2}(p)-p\right\| &=&\left\|
  I_{2}p-I_{2}a_{2}+a_{2}-p\right\|  \\ &\leq&
2\epsilon
\end{eqnarray*}
and so
\begin{eqnarray*}
\left\| E_{2}(a_{2})-p\right\| &=&\left\|
E_{2}(a_{2})-E_{2}(p)+E_{2}(p)-p\right\|  \\&\leq& \left\| E_{2}\right\|
\left\|
a_{2}-p\right\| +\left\| E_{2}(p)-p\right\| \\& \leq& 3\epsilon \\&<&1,
\end{eqnarray*}
and since $E_{2}(a_{2})\in \mathcal{C}_{2}$, we have a contradiction.
\end{proof}
\begin{lma}
Let $\mathcal{A}$ be an operator algebra with masa $\mathcal{C}.$
Suppose there exists a chain
of subalgebras with masas ($\mathcal{A}_{1},\mathcal{C}_{1})\subseteq
$ ($\mathcal{A}_{2},\mathcal{C}_{2})\subseteq\ldots$ whose union is
dense in $\mathcal{A}$. Then $\mathcal{C}=cl(
\bigcup _{k}\mathcal{C}_{k}).$
\end{lma}
\begin{proof}
Let $\left\{ e_{i,j}^{k}\right\} _{i,j}$ be a matrix unit system for
$\mathcal{A}_{k}$
compatible with the masa $\mathcal{C}_{k.}$ Take any $\epsilon >0$ and
$c\in \mathcal{C}.$ Then
we can select a sufficiently large index, $k$ and an element $a_{k}$ of
$%
\mathcal{A}_{k}$ with $\left\| c-a_{k}\right\| <\epsilon .$ Then
$\left\| E_{k}(c)-E_{k}(a_{k})\right\| <\epsilon $ and $\left\|
c-E_{k}(c)\right\| <2\epsilon ,$ thus $\left\| c-E_{k}(a_{k})\right\|
<3\epsilon $ and since $E_{k}(a_{k})\in \mathcal{C}_{k},$ we are done.
\end{proof}
\begin{lma}
Let
$(\mathcal{A}_{1},\mathcal{C}_{1}),(\mathcal{A}_{2},\mathcal{C}_{2})$ be
digraph subalgebras of $(\mathcal{A},\mathcal{C})
$, such that $N_{\mathcal{C}_{i}}(\mathcal{A}_{i})\subseteq
N_{\mathcal{C}}(\mathcal{A})$, $i=1,2.$ Then given $%
\delta >0$ there exists $\epsilon >0$ such that if $\mathcal{A}_{1}
\subseteq_{\epsilon} \mathcal{A}_{2},$ then
$N_{\mathcal{C}_{1}}(\mathcal{A}_{1})\subseteq_{\delta} %
N_{\mathcal{C}_{2}}(\mathcal{A}_{2}).$
\end{lma}
\begin{proof}
Take $v\in N_{\mathcal{C}_{1}}(\mathcal{A}_{1}).$ Then $v\in
N_{\mathcal{C}}(\mathcal{A})$ and thus for all $c\in \mathcal{C}$
with $\left\| c\right\| \leq 1$ there exists $d\in \mathcal{C}$ with
$\left\| v^{\ast
}cv-d\right\| =0$ and similarly for $vcv^{\ast }.$ Since
$\mathcal{A}_{1}\subseteq_{\epsilon}\mathcal{A}_{2}$ we can find $w\in
ball(\mathcal{A}_{2})$ with $\left\|
v-w\right\| <\epsilon .$ Our aim is to show that $w$ is an approximate
partial isometry which approximately normalises $\mathcal{C}_2$.  The
former property
is clear from the fact that $w$ is close to $v$, an exact partial
isometry.
Also $w$ approximately normalises $\mathcal{C}$, again because it is
close to the
exactly normalising element $v$, so for all $c\in ball(\mathcal{C})$
there exist $%
c_{1},c_{2}\in \mathcal{C}$ with $\left\| w^{\ast }cw-c_{1}\right\|
<\delta
_{1}(\epsilon ),$ $\left\| wcw^{\ast }-c_{1}\right\| <\delta
_{1}(\epsilon ).
$ Now take $c_{2}\in ball(\mathcal{C}_{2}).$ Since
$\mathcal{C}_{2}\subseteq \mathcal{C},$ $v^{\ast
}c_{2}v\in \mathcal{C}.$ We first need to estimate $\left\| v^{\ast
}c_{2}v-E_{2}(v^{\ast }c_{2}v)\right\| .$ Now
\begin{eqnarray*}
E_{2}(v^{\ast }c_{2}v)&=&\underset{i}{\sum }e_{i}^{2}v^{\ast
  }c_{2}ve_{i}^{2} \\&=&
\underset{i}{\sum }e_{i}^{2}v^{\ast }c_{2}v \\&=&I_{2}v^{\ast }c_{2}v,
\end{eqnarray*}
\ since  $e_{i}^{2},v^{\ast
}c_{2}v\in \mathcal{C}$
and where $I_{2}=\underset{i}{\sum }e_{i}^{2}.$ Then
\begin{eqnarray*}
\left\| v^{\ast }c_{2}v-E_{2}(v^{\ast }c_{2}v)\right\| &=&\left\|
v^{\ast
}c_{2}v-I_{2}v^{\ast }c_{2}v\right\|  \\&=&\left\| v^{\ast
}c_{2}v-w^{\ast
}c_{2}w+I_{2}w^{\ast }c_{2}w-I_{2}v^{\ast }c_{2}v\right\| \\&\leq&
2\delta
_{1}(\epsilon ).
\end{eqnarray*}
Now $E_{2}(w^{\ast }c_{2}w)\in \mathcal{C}_{2}$ and we estimate
\begin{eqnarray*}
\left\| w^{\ast }c_{2}w-E_{2}(w^{\ast }c_{2}w)\right\| &=&\left\|
w^{\ast
}c_{2}w-v^{\ast }c_{2}v+v^{\ast }c_{2}v-E_{2}(v^{\ast
}c_{2}v)+E_{2}(v^{\ast
}c_{2}v)-E_{2}(w^{\ast }c_{2}w)\right\| \\
&\leq&4\delta_{1}(\epsilon).
\end{eqnarray*}
Thus $w$ approximately normalises $\mathcal{C}_{2}.$ Lemma 3.5 now
provides a
normalising partial isometry $\hat{w}\in
N_{\mathcal{C}_{2}}(\mathcal{A}_{2})$ with $\left\|
w-\hat{w}\right\| <4\delta _{1}(\epsilon ).$ Thus
$N_{\mathcal{C}_{1}}(\mathcal{A}_{1})%
\subseteq_{4\delta _{1}(\epsilon
)}N_{\mathcal{C}_{2}}(\mathcal{A}_{2}).$
\end{proof}
\begin{lma}
Let$\;($$\mathcal{A}_{1},\mathcal{C}_{1}),(\mathcal{A}_{2},\mathcal{C}_{2})$

be digraph algebras. Let $0<\epsilon <%
\frac{1}{4}$ be given and suppose $\mathcal{C}_{1}\subseteq
\mathcal{C}_{2}.$ Then there exists $%
\delta >0$ dependent only on $\epsilon $ and $\mathcal{A}_{1}$ such that
if $%
N_{\mathcal{C}_{1}}(\mathcal{A}_{1})\subseteq_{\delta}N_{\mathcal{C}_{2}}(\mathcal{A}_{2}),$
 then there
exists a regular star extendible algebra injection $\phi
:\mathcal{A}_{1}\rightarrow
\mathcal{A}_{2}$ with $\left\| id-\phi \right\| <\epsilon .$
\end{lma}
\begin{proof}
Firstly, we observe that if $v_{1},v_{2}\in M_{n}$ are permutation type
partial isometries with $\left\| v_{1}-v_{2}\right\| <1,$ then $v_{1}$
and $%
v_{2}$ have the same support, for if not we could find minimal
projections $%
p\;$and $q$ with $1=\left\| p(v_{1}-v_{2})q\right\| \leq \left\|
v_{1}-v_{2}\right\| <1.$ Secondly if $v_{2}v_{2}^{\ast }=v_{1}^{\ast
}v_{1},$
then $v_{1}v_{2}$ is another permutation type partial isomerty. With
these
preliminary observations in mind, we proceed with the proof. Suppose
$\mathcal{A}_{1}$
is a digraph algebra on $n_{1}$ vertices, then any cycle within the
digraph
will have length no greater than $n_{1}.$ Let $G_{s}$ be any spanning
tree
for the digraph of $\mathcal{A}_{1}$ fixed throughout, and let $%
e_{1},e_{2},...,e_{n_{1}-1}$ be the matrix units corresponding to the
edges
of $G_{s}.$ Since $G_{s}$ has no cycles, any matrix unit in
$\mathcal{A}_{1}$ can be
written uniquely as a word of minimal length, no greater than $n_{1},$
using the
alphabet $E=\left\{ e_{1},e_{2},...,e_{n_{1}-1},e_{1}^{\ast
},e_{2}^{\ast
},...,e_{n_{1}-1}^{\ast }\right\} .$ To fix notation, let $\left\{
e_{i,j}\right\} _{i,j},\left\{ g_{k,l}\right\} _{k,l}$ be matrix unit
systems for $\mathcal{A}_{1}$ and $\mathcal{A}_{2}$ respectively,
compatible with the given
masas. Since $\mathcal{C}_{1}\subseteq \mathcal{C}_{2}$, each diagonal
matrix unit $e_{i,i}$ can
be written as a sum of the $g_{l,l}$'s. Denote this sum for
each $e_{i,i}$ by $f_{i,i},$ for $1\leq i\leq n_{1}.$ Now take the
matrix
unit $e_{1}.$ Since $e_{1}\in N_{\mathcal{C}}(\mathcal{A})$ there exists
$v_{1}\in
N_{\mathcal{C}_{2}}(\mathcal{A}_{2})$ with $\left\| e_{1}-v_{1}\right\|
<\delta .$ If $%
e_{1}^{\ast }e_{1}=e_{i_{1},i_{1}},$ and  $e_{1}e_{1}^{\ast
}=e_{j_{1},j_{1}}$
then, $\left\| v_{1}^{\ast }v_{1}-f_{i_{1},i_{1}}\right\| <2\delta $ and
$%
\left\| v_{1}v_{1}^{\ast }-f_{j_{1},j_{1}}\right\| <2\delta .$ Since $%
f_{i_{1},i_{1}}$ and $f_{j_{1},j_{1}}$ are both standard projections and
$%
v_{1}$ is a unimodular sum of the $g_{k,l}$'s , provided $\delta
<\frac{1}{2},$ $v_{1}v_{1}^{\ast }=f_{j_{1},j_{1}}$ and $v_{1}^{\ast
}v_{1}=f_{i_{1},i_{1}}$ and so $v_{1}$ has the correct initial and final
projections. Now set $f_{1}=v_{1}.$ Similarly, we create $%
f_{2},...,f_{n_{1}-1}$ having the right initial and final projections
(in
the above sense) and with $\left\| e_{i}-f_{i}\right\| <\delta $ for
each $i.
$ We now form the corresponding alphabet $F=\left\{
f_{1},f_{2},...,f_{n_{1}-1},f_{1}^{\ast },f_{2}^{\ast
},...,f_{n_{1}-1}^{\ast }\right\} .$ Now take any other matrix unit
$e_{i,j}$
in $A_{1}$ and let $w_{E}^{i,j}$ denote its unique word of minimal
length in $E.$ Define $%
f_{i,j}$ to be the element with corresponding word in $F.$ Note that $%
\left\| e_{i,j}-f_{i,j}\right\| =\left\| w_{E}^{i,j}-w_{F}^{i,j}\right\|
\leq n_{1}\delta .$ We now need to show that $f_{i,j}\in
\mathcal{A}_{2}$ for each $%
i,j.$ Firstly by construction, each $f_{i,j}$ is a permutation type
partial
isometry in the containing $C^{\ast }$-algebra. Since $e_{i,j}\in
N_{\mathcal{C}_{1}}(\mathcal{A}_{1})$ we can find $v_{i,j}\in
N_{\mathcal{C}_{2}}(\mathcal{A}_{2})$ with $\left\|
e_{i,j}-v_{i,j}\right\| <\delta .$ Then 
\begin{eqnarray*}
\left\| f_{i,j}-v_{i,j}\right\|
&=&\left\| w_{F}^{i,j}-v_{i,j}\right\| \\&\leq& \left\|
w_{F}^{i,j}-w_{E}^{i,j}\right\| +\left\| e_{i,j}-v_{i,j}\right\|
 \\&<&n_{1}\delta +\delta .
\end{eqnarray*}
Provided we choose $\delta $ sufficiently small so that $n_{1}\delta
+\delta
<1$, $f_{i,j}$ must have the same support as $v_{i,j},$ and thus
$f_{i,j}\in
\mathcal{A}_{2}.$ In this way we create a matrix unit system $\left\{
f_{i,j}\right\}
_{i,j}$ for a copy of $\mathcal{A}_{1}$ in $\mathcal{A}_{2}$ with
$\left\|
e_{i,j}-f_{i,j}\right\| \leq n_{1}\delta \;$for each $i,j.$ Set $\phi $
to be the linear extension of the
correspondences $e_{i,j}\rightarrow f_{i,j}.$ Then $\phi $ is 
regular, since $\phi (N_{\mathcal{C}_{1}}(\mathcal{A}_{1}))\subseteq
N_{\mathcal{C}_{2}}(\mathcal{A}_{2})$ and it is clear that given
$\mathcal{A}_{1}$ we can choose $\delta$ sufficiently small so that
$\left\| id-\phi \right\| <\epsilon .$
\end{proof}
\begin{lma}
Let
$(\mathcal{A}_{1},\mathcal{C}_{1}),\;(\mathcal{A}_{2},\mathcal{C}_{2})$
be digraph subalgebras of an operator
algebra $(\mathcal{A},\mathcal{C})$ such that
$N_{\mathcal{C}_{i}}(\mathcal{A}_{i})\subseteq
N_{\mathcal{C}}(\mathcal{A})$, $i=1,2.$
Then, given $\epsilon >0$ we can find $\delta >0$ such that if
$\mathcal{A}_{1}%
\subseteq_{\delta}\mathcal{A}_{2},$ there exists a regular star
extendible
algebra injection $\phi :\mathcal{A}_{1}\rightarrow \mathcal{A}_{2}$
with $\left\| id-\phi
\right\| <\epsilon ,$ where $\delta $ depends only on $\mathcal{A}_{1}$
and $\epsilon .
$
\end{lma}
This lemma may be viewed as the regular analogue of stability. Put
more succinctly it says that the family of digraph algebras is regularly
stable.
\begin{proof}
By Lemma 3.6, we choose $\delta <\frac{1}{4}$ so that
$\mathcal{C}_{1}\subseteq \mathcal{C}_{2}.$
Now Lemma 3.9 implies the existence of $\delta _{1}>0$ such that if
$\mathcal{A}_{1}%
\subseteq_{\delta _{1}}\mathcal{A}_{2}$ and
$N_{\mathcal{C}_{1}}(\mathcal{A}_{1})\subseteq_{
\delta _{1}} N_{\mathcal{C}_{2}}(\mathcal{A}_{2})$ there exists a
regular star
extendible injection $\phi :\mathcal{A}_{1}\rightarrow \mathcal{A}_{2}$
with $\left\| id-\phi
\right\| <\epsilon .$ By Lemma 3.8 there exists $\delta _{2}$ for which
given $%
\mathcal{A}_{1}\subseteq_{\delta _{2}} \mathcal{A}_{2}$ we have
$N_{\mathcal{C}_{1}}(\mathcal{A}_{1})%
\subseteq_{\delta _{1}} N_{\mathcal{C}_{2}}(\mathcal{A}_{2}).$
\end{proof}
We have now arrived at the promised characterisation of regular limits
of digraph algebras.
\begin{thm}
A separable operator algebra with masa,
$(\mathcal{A},\mathcal{C}_{\mathcal{A}})$ is a regular limit of
digraph  algebras if and only if for each $\epsilon >0$ and finite
subset $%
\mathcal{S}\subseteq \mathcal{A}$, there exists a pair
$((\mathcal{B},\mathcal{C}_{\mathcal{B}}),\phi )$ of digraph algebra
and regular star extendible injection $\phi :\mathcal{B}\rightarrow
\mathcal{A}$ with $%
dist(\mathcal{S},\phi (\mathcal{B}))<\epsilon .$
\end{thm}

\begin{proof}
Necessity of the local condition is clear, so suppose the local
condition
holds. Choose a dense sequence in the unit ball of $\mathcal{A},$
$\left\{
a_{k}\right\} .$ Let $\left\{ \epsilon _{k}\right\} $ be a summable
sequence, with $\epsilon _{k}<\frac{1}{4}$ for each $k.$ By hypothesis,
there exists a pair $((\mathcal{A}_{1},\mathcal{C}_{1}),\phi _{1})$ with
$\phi
_{1}(N_{\mathcal{C}_{1}}(\mathcal{A}_{1}))\subseteq
N_{\mathcal{C}_{\mathcal{A}}}(\mathcal{A})$ and $dist(a_{0},\phi
_{1}(\mathcal{A}_{1}))<\epsilon _{0}.$ Now, given $\phi
_{1}(\mathcal{A}_{1})$ and $\epsilon _{1}
$, Lemma 3.10 implies that we can find $\delta _{1}>0$ such that if $%
(\mathcal{A}_{2},\mathcal{C}_{2})$ is another digraph algebra such that
$N_{\mathcal{C}_{2}}(\mathcal{A}_{2})%
\subseteq N_{\mathcal{C}_{\mathcal{A}}}(\mathcal{A})$ and $\phi
_{1}(\mathcal{A}_{1})\subseteq_{\delta _{1}}
\mathcal{A}_{2},$ then there exists a regular star extendible algebra
injection $\pi
_{1}:\mathcal{A}_{1}\rightarrow \mathcal{A}_{2}$ with $\left\| id-\pi
_{1}\right\| <\epsilon
_{1}.$ We now demonstrate how the local condition provides
$\mathcal{A}_{2}.$ Since $%
\mathcal{A}_{1}$ is finite dimensional, we can select a finite
$\frac{\delta _{1}}{2}$
net for the unit ball of 
$\phi _{1}(\mathcal{A}_{1}),\;\left\{
a_{1}^{1},...,a_{r_{1}}^{1}\right\} ,$we assume $\delta _{1}<\epsilon
_{1}$.
Consider the finite subset 
$S=\left\{a_{0},a_{1},a_{1}^{1},...,a_{r_{1}}^{1}\right\} $ $\subseteq
\mathcal{A}.$ By the
local condition there exists a digraph algebra and a regular star
extendible
injection $\phi _{2}:\mathcal{A}_{2}\rightarrow \mathcal{A}$ with
$dist(a_{i},\phi
_{2}(\mathcal{A}_{2}))<\delta _{1}$ $i=0,1$ and $dist(a_{i}^{1},\phi
_{2}(\mathcal{A}_{2}))<%
\frac{\delta _{1}}{2}$ for $i=1,...,r_{1},$ from which $\phi
_{1}(\mathcal{A}_{1})%
\subseteq_{\delta _{1}} \phi _{2}(\mathcal{A}_{2})$ and
$\mathcal{C}_{1}\subseteq
\mathcal{C}_{2}.$ Continuing in this way we construct a sequence of
finite dimensional operator algebras $%
(\mathcal{A}_{k},\mathcal{C}_{k})$ and regular star extendible
injections $\pi
_{k}:\mathcal{A}_{k}\rightarrow \mathcal{A}_{k+1}$ with $\left\| id-\pi
_{k}\right\| <\epsilon
_{k}$ and $\mathcal{C}_{k}\subseteq \mathcal{C}_{k+1}$ for each $k.$ Now
consider the
diagram
\begin{equation*}
\begin{CD}
\mathcal{A}_1 @>\pi_1>> \mathcal{A}_2 @>\pi_2>> \mathcal{A}_3 @>\pi_3>>
\mathcal{A}_4 @>\pi_4>> \cdots  \\
 \symbdown{I_1} & & \symbdown{I_2} &  & \symbdown{I_3}
      &  & \symbdown{I_4} & 
      & & \\
\mathcal{A} @>id>> \mathcal{A} @>id>> \mathcal{A} @>id>> \mathcal{A}
@>id>>  \cdots & 
\end{CD} 
\end{equation*}
where each $I_{k}$ is the restriction to $\phi _{k}(\mathcal{A}_{k})$ of
the identity
map. We estimate
\begin{equation*}
\left\| I_{k+1}\circ \pi _{k}-id\circ I_{k}\right\| <\epsilon _{k} 
\end{equation*}
for
each $k.$
Since $\left\{ \epsilon _{k}\right\} $ is summable, the diagram commutes
asymptotically, thus $\mathcal{A}=\lim (\phi _{k}(\mathcal{A}_{k}),\pi
_{k})$ and $\mathcal{C}=\lim (\mathcal{C}_{k})
$ and the proof is complete.
\end{proof}
\textit{Acknowledgements} The author is supported by an EPSRC
studentship. I would like to thank my supervisor Prof S. C. Power for
his invaluable advice and encouragement while I was writing this paper,
Dr. R. M.  Green for his many useful suggestions and Dr. J. S. Groves.

\vspace{0.5cm}Department of Mathematics and Statistics\\
Lancaster University\\
Lancaster\\
LA1 4YF\\
UK\\
email: p.haworth@lancaster.ac.uk\\
Fax: +44 01524 592681

\end{document}